\documentstyle[12pt]{article}
\begin{document}

\title{On the Metric Independent Exotic Homology.}

\author{S.Novikov}

\date{
\centerline{IPST, University of Maryland, College Park,}
\centerline{MD,20742-2431, USA}
\centerline{novikov@ipst.umd.edu} 
}

\maketitle

\begin{abstract}
Different types of nonstandard homology groups based on the
various subcomplexes  of differential forms are considered as a
continuation of the recent authors works. Some of them reflect
interesting properties of dynamical systems on the compact
manifolds. In order to study them a Special Perturbation Theory in
the form of Spectral Sequences is developed. In some cases a
convenient fermionic formalism of dealing with differential forms
is used originated from the work of Witten in the Morse Theory
(1982)and the authors work  where some nonstandard analog of Morse
Inequalities for vector fields was found (1986).
\end{abstract}

In the work \cite{N3} we invented some sort of exotic homology of
the first and second kind. Let us remind here that the second kind
exotic homology defined in this work are constructed in  the
following natural way: for every linear space $L$ with operator
$d':L\rightarrow L$ we can define homology on any $d$-invariant
subspace $T\subset L$ such that $(d')^2:T\rightarrow 0$. So we
have by definition $H_T=Ker(d')/Im(d')$. Our main example was
based on the De-Rham complex $\Lambda^*(M)=\sum_{i=0}^n\Lambda^i$
of real-valued differential forms for any $C^{\infty}$-manifold
$M$. We considered a family of operators $d'=d+\lambda\omega^*$
such that $d'(a)=da+\lambda\omega\wedge a$. Here $\omega$ is an
one-form (maybe non closed), and $a$ is any $C^{\infty}$-form. A
lot of work was done  since 1986  (see, for example,  \cite{N1,P})
for the case of the closed one-form $\omega$ where $(d')^2=0$. We
do not discuss this case here.

I.{\it Perturbation of $d$ by the nonclosed odd form.}

 Consider the standard real-valued nonclosed forms
$\omega$.  Our operators $d'$ are defined for all values of
$\lambda\in R$ on the subspace $T\subset \Lambda^*$ such that
$a\in T$ if and only if $(\Omega)\wedge a=(\Omega)^*(a)=0$ where
$\Omega=d\omega$. Here $R$ is a field of ordinary numbers like $R$
or $C$. We especially mention that here because later on (see the
part II) this ring will became a $Z_2$-graded supercommutative
ring with unit. Let us concentrate first on the special case when
$\lambda=0$ but the standard De Rham operator $d'=d$ is considered
on the more general class of subspaces $T\subset \Lambda^*(M)$
such that $a\in T$ if and only if $a\wedge\Omega=0$. Here $\Omega$
is any fixed closed differential form $d\Omega=0$. In the example
above we have $\Omega=d\omega$, so this form is exact. Anyway, the
operator $d$ commutes with the multiplication operator
$$\Omega:u\rightarrow u\wedge\Omega, d\Omega=\Omega d$$ We take
now an odd-dimensional (2n+1)-manifold $ M$ which is a
nondegenerate energy level in the symplectic manifold $N$ with
symplectic 2-form $w$. We take a $2n$-form $\Omega=w^n$ restricted
on the $2n+1$-submanifold $M$. It turns out that the exotic
homology here are associated with Hamiltonian System (1-foliation)
on the submanifold $M$. We denote the subspace $T$ here as
$T_{\Omega}$ and homology $Ker(d)/Im(d)$ in this subspace as
$H^*_{\Omega}$. The form $w$ restricted on $M$ is degenerate along
the Hamiltonian vector field $X$ only, i.e. $w(X,Z)=0$ for every
vector field $Z$. People in Symplectic Geometry call $X$ a ''Reeb
Vector Field''.

\newtheorem{lem}{Lemma}
\begin{lem} $C^{\infty}$-differential k-form $a\in \Lambda^*(M)$
 belongs to the subspace $T_{\Omega}^k\subset\Lambda^k$ if and only if

a. It is equal to zero for  $k=0$;

 b. For $k=1$ its values on the
vector field $X$ is identically equal to zero;

c. For $k\geq 2$ the subspace $T_{\Omega}^k$ is equal to the whole
space $\Lambda^k(M)$
\end{lem}

This statement was proved in the work \cite{N3} for $n=1$ in the
slightly different terminology. For the case $k=0$ our statement
is obvious. There is no difference for the case $k=1,n>1$. Let us
remind this proof here. According to the Darboux theorem, there
exist local coordinate system $p_i,q_i,r, i=1,2,\ldots,n$ such
that $w=\sum_idp_i\wedge dq_i$. So the direction of the Reeb
vector field  exactly corresponds to the variable $r$. We have
$w^n=(const)\prod dp_1\wedge\ldots\wedge dq_n$, so the equality
$w^n\wedge a=0$ implies locally $a=\sum_i (c_idp_i+l_idq_i)$
without $dr$. Therefore we have $(a,X)=0$ because
$(dp_i,X)=(dq_i,X)=0$. For $k=1$ our lemma is proved. For $k>1$ it
is obvious.  Consider now the homology of this complex
$$0\rightarrow
T^1_{\Omega}\rightarrow\Lambda^2(M)\rightarrow\Lambda^3(M)\rightarrow\ldots\Lambda^n(M)
\rightarrow 0$$ with standard boundary operator $d$.

\newtheorem{th}{Theorem}
\begin{th} For the compact smooth manifold $M$  a natural exact
sequence is well-defined:

$$0\rightarrow R\rightarrow Ker(\nabla_X)\rightarrow
H^1_{\Omega}\rightarrow H^1(M)\rightarrow
C^{\infty}(M)/\nabla_X(C^{\infty}(M))\rightarrow
H^2_{\Omega}\rightarrow H^2(M)\rightarrow 0$$

where $$\nabla_X:C^{\infty}(M)\rightarrow C^{\infty}(M)$$ is a
derivative of function $f$ along vector field $X$
$$\nabla_X(f)=X^i\partial_if, X=(X^i)$$

\end{th}

{\it Remark: I clarified this question after the very useful
discussion with D.Dolgopyat.}

 Let us construct all homomorphisms
and prove exactness of this sequence for compact manifolds. We
start from the left part.

For every function $f\in C^{\infty}(M)$ we have a form $df$. The
equality $\nabla_Xf=X^i\partial_if=0$ implies $(df,X)=0$ by
definition. So we have a map $Ker(\nabla_X)\rightarrow
H^1_{\Omega}$ whose kernel is exactly constant functions.  The
image of this map in the group $H^1_{\Omega}$ consists of all
exact forms such that $(u,X)=0$. Therefore a factor-group by this
image lies in the first homology group $H^1(M)$ Therefore we
constructed a second map $H^1_{\Omega}\rightarrow H^1(M)$ and
proved exactness of the sequence in the term $H^1_{\Omega}$.

For the construction of the next map, we simply take $u\rightarrow
(u,X)$. This function belongs to $C^{\infty}(M)$. Varying the
closed 1-form $u$  in the homology class $u+df$ we have a correct
result in the factor $C^{\infty}(M)/\nabla_X(C^{\infty}(M))$. Its
kernel contains exactly all closed 1-forms $u$ such that there
exists  a function $f$ for which we have $(u+df,X)=0$. So we
proved exactness in this term.

By definition, the next group $H^2_{\Omega}$ consists of all
closed 2-forms $dv=0$ modulo $du$ where $(u,X)=0$. Let us consider
part $H^{2,0}_{\Omega}\subset H^2_{\Omega}$ of this group
represented by the exact 2-forms $v=dz$ but maybe $(z,X)\neq 0$.
The set of the quantities $(z,X)\in C^{\infty}(M)$ can be
identified with all space $C^{\infty}(M)$. So we have a map
$$C^{\infty}(M)\rightarrow H^{2,0}_{\Omega}\subset H^2_{\Omega}$$
whose kernel is exactly represented by the  projections of the
closed 2-forms $(u,X)\in C^{\infty}(M)$. This projection is
well-defined modulo terms like $\nabla_Xf$. So this map is also
constructed. Our sequence is exact in this term as well.

The next map $H^2_{\Omega}\rightarrow H^2(M)$ is natural. Its
kernel obviously is equal to $H^{2,0}_{\Omega}$ by definition. It
is an epimorphism.

So our theorem is proved. By definition, for the case $\lambda=0$
we have $$H^i_{\Omega}=H^i(M), i\geq 3$$

  {\it Question: We
defined our exact sequence formally for smooth functions only. Is
it possible to make a proper complition of it to the Hilbert
spaces such that all homology will remain the same as in the Hodge
Theory? Our complex is nonelliptic here. Is it true that after the
proper ergodicity requirements our homology will be
finite-dimensional even after the complition?}

Let us consider now the case $\lambda \neq 0$ where
$\Omega=d\omega$. We are dealing with the same subspace in the
space of smooth differential forms, but $Z$-grading is lost: it
should be replaced by the $Z_2$-grading
$H^{odd}_{\Omega,\lambda\omega}$ and
$H^{even}_{\Omega,\lambda\omega}$ where
$H^*_{\Omega,\lambda=0}=H^*_{\Omega}$. Acting by the operator
$d'=d+\lambda \omega^*$ on the odd forms
$v_1+v_3+v_5+\ldots=v+\lambda
v'+\lambda^2v''+\ldots\in\Lambda^{odd}(M)$ and even forms
$u_2+u_4+\ldots=u+\lambda u'+\lambda^2
u''+\ldots\in\Lambda^{even}(M)$, we are coming along the line of
the work \cite{N1} to the spectral sequence:  its first
differential is determined by the zero order in $\lambda$, so it
is ordinary $d$. Its  second differential (determined by the first
order approximation in the variable $\lambda$) is generated by the
multiplication on the form $\omega$. In this simple case this
spectral sequence has only two nontrivial differentials $d_1=d,
d_2=\omega^*$ for all dimensions $2n+1>3$:

\begin{lem} Following differential is well-defined on the group
$H^*_{\Omega}=E_2$:
$$d_2=\omega^*:H^{odd}_{\Omega,\lambda=0}\rightarrow
H^{even}_{\Omega,\lambda=0}\rightarrow H^{odd}_{\Omega,\lambda=0},
d_2^2=0 $$ such that $$d_2(u)=\omega^*(u)=\omega\wedge u$$ for the
even forms, and $$d_2(v)=\omega^*(v)=\omega\wedge v$$ for the odd
forms,  where $du=dv=0$ are the even and odd representative
cocycles correspondingly.
\end{lem}

{\it Remark: The groups $Ker(\omega^*)/Im(\omega^*)$ in good cases
are equal to $H^*_{\Omega,\lambda\omega}$ for odd and even cases
for all odd dimensions $n>3,\lambda\neq 0$ and small enough.
However, in order to prove this statement we need to perform  more
serious analysis of this complex taking into account its
nonellipticity.}

For the special case $n=3$ our spectral sequence is more
complicated: we define ''Massey Products'' generated by the form
$\omega$: By definition, we call by the first Massey product a
multiplication operator $\omega^*:u\rightarrow \omega\wedge u$.
The second Massey product is a following operator: let
$\omega\wedge u=dv$ in our complex. We put
$$d_3(u)=\{\omega,\omega,u\}=\omega\wedge v$$ This operation is
well-defined on the homology group $E_3$ determined by the
operator $d_2=\omega^*$ on the group $E_2=H^*_{\Omega,\lambda=0}$
with natural $Z_2$-grading. By induction, we define the next
Massey product $d_{l+1}$ acting on the homology group
$E_{l+1}=H(E_l,d_l)$ of the operator $d_l$ acting on $E_l$ with
natural $Z_2$-grading as a map reversing  this grading by the
formula: $$d_{l+1}(u)=\{\omega,\ldots,\omega,u\}=\omega\wedge v$$
( here $\omega$ enters $l$ times in the $l$-th Massey Product
equal to the differential $d_{l+1}$ ), such that $dv=d_l(u)$, and
all $d_j(u)$ are equivalent to $0$ for $j\leq l$. Let us mention
that we are speaking here about the class of $u$ in the group
$E_l$.

 So we consider any compact odd-dimensional manifold $M=M^{2n+1}$ in the
symplectic (noncompact)  manifold $N=N^{2n+2}$ with exact
symplectic 2-form $d\omega$, and operator $d'=d+\lambda\omega$ on
$M$ acting in the kernel of the multiplication operator by the
form $\Omega=(d\omega)^n$.

\begin{th} For $2n+1=3$ the sequence on Massey products $\{\omega,\ldots,\omega,u\}$
determines a well-defined spectral sequence of groups and
differentials $E_l=E^{odd}\bigoplus E^{even}_l,d_l$ where
differentials  are grading reversing. The group $E_1$ are equal to
the kernel of multiplication operator by the form $\Omega=d\omega$
in the space of all smooth differential forms $\Lambda^*(M)$, and
$d_1=d$. For all odd dimensions $2n+1>3$ all differentials $d_i$
are equal to zero for $i>2$.
\end{th}

Of course, this Spectral Sequence is a partial case of the more
general class of spectral sequences generated by the pair of two
anticommuting differentials $d'=d_1+d_2,
(d')^2=d_1^2+d_2^2+d_1d_2+d_2d_1=0$, but in our cases all spectral
sequence can be expressed in the form of Massey Products. This
analytically defined spectral sequence and its calculation through
the Massey Products appeared first time in the work \cite{N1},
 for the closed 1-form $\omega$ where we were dealing with the
elliptic complex of all forms. As M.Farber pointed out to me, this
construction solved the problem of calculation of differentials in
the purely homological spectral sequence already invented by him
purely topologically in the work \cite{F2}, developing an old work
of J.Milnor \cite{M}, so our analytical approach very easily led
to the results unclear from the purely topological point of view.
 The general existence theorem of the analytical ''perturbation
spectral sequence'', was formulated and proved for the
deformations of all elliptic complexes in \cite{F}. We studied in
details some specific very interesting example in the work
\cite{N3}. In the present work we are dealing with nonelliptic
complexes, so we may have difficulties of the functional type. We
cannot identify at the moment the groups $E_{\infty}$ for this
spectral sequence. As B.Mityagin pointed out to the author, John
von Neumann  in 1930s  performed some considerations looking like
construction of our perturbation spectral sequence in the theory
of operators, but there was no notions of the homological algebra
at that time for the proper understanding what is going on. This
comparison was done first time in the work \cite{N3}.

{\bf More general exotic homology}
 can be naturally defined for the differential forms $\Lambda^*(M,R)$
with
coefficients in the $Z_2$-graded skew commutative
(supercommutative) associative ring with unit $R=R^+\bigoplus R^-$
instead of ordinary numbers: they are based on the operator
$$d'=d+\sum_{i\geq 1}\lambda_i\omega^{i*}$$ where dimension of the
form $\omega^i$ is equal to $i\geq 1$, and $\lambda_i\in R^+$ if
$i=2j+1, \lambda_i\in R^-$ if $i=2j$. We have
$(\sum_i\lambda_i\omega^i)^2=0$.

 These exotic homology groups are
defined on the subspace $T_{\Omega}\subset \Lambda^*(M)$ where
$\Omega=\sum_i\lambda_i(d\omega^i)^*$. We have
$$T_{\Omega}=T^{even}_{\Omega}\bigoplus T^{odd}_{\Omega}$$
$$d':T^{even}\rightarrow T^{odd}\rightarrow T^{even}$$ So we have
finally $H^{odd}_{\Omega}$ and $H^{even}_{\Omega}$. Many previous
statements can be easily extended to this case.

\vspace{0.3cm}

II.{\it Perturbation of $d$ by the contraction operators}.

We can define another class of exotic homology based on the
perturbation of the previously defined operators  by the following
tensor contraction operators on the spaces of differential forms
with coefficients in the supercommutative $Z_2$-graded ring $R$:
consider the operator $d'$ $$d'=d+\lambda_i\sum_{i\geq
1}\omega^{i*} + \sum_j \mu_j\hat{X}_j$$ where the last operators
$\hat{X}_j$ are equal to the action of some selected skew
symmetric tensors $X_j$ with $j$ upper indices on differential
forms by the purely algebraic standard tensor contraction. By
definition, the action of operator $\hat{X}_j$ on the forms of
dimension less than $j$ is equal to zero, coefficients
$\lambda_i\in R^+$ belong to the even part of the ring $R$ and
$\mu\in R^-$ belong to the odd part of $R=R^+\bigoplus R^-$. The
simplest case here is an operator
$$d'=d+\lambda\omega^{1*}+\mu\hat{X}_1$$ where $\lambda,\mu$ are
the ordinary numbers. Let us remind here that there exists a
natural external product of skew symmetric tensors with upper
indices as well as for differential forms. Following
\cite{W,NS,S,N2}, we use an algebraic language of fermionic
creation and annihilation operators very convenient for the
calculations with differential forms and other skew-symmetric
tensors. In the given system of local coordinates $x^i$ we
introduce operators $a^i,a^+_j$ with standard commutation
relations $$a^ia^+_j+a^+_ja^i=\delta^i_j,a^+_ia^+_j=-a^+_ja^+_i,
a^ja^i=-a^ia^j$$ No metric on the manifold is given, so we do not
consider annihilation operators as conjugated to the creation
operators--no conjugation operation can be defined without metric.
The local basis of differential $k$-forms is given in this
''Dirac-Fock Space'' by the creation operators applied to the
''vacuum vector''
 $$dx^{i_1}\wedge\ldots \wedge dx^{i_k}=a^{i_1}a^{i_2}\ldots a^{i_k}\Phi_0$$
''Vacuum''  vector  represents a
 constant function $1$, i.e. ($\Phi_0=1$). It
satisfies to the relations $a^+_j\Phi_0=0$. The operators
$\hat{X}^{i_1\dots i_k}$ corresponding to the contraction with
skew symmetric tensor fields $X^{i_1\ldots i_k}(x)$ with $k$ upper
indices are given by the operators $$\hat{X}^{i_1\ldots
i_k}=X^{i_1\ldots i_k}(x)a^+_{i_1}\ldots a^+_{i_k}$$ acting in
this Fock space. Let $\omega^1=\omega=\sum_i\omega_idx^i$ and
$X_1=X=(X^i(x)$. Denoting $\omega^1$  by $\omega$ and $X_1$ by
$X$, we have:
$$d=a^i\partial_i,\omega^*=\omega_i(x)a^i,\hat{X}_1=\hat{X}=X^i(x)a^+_i$$
We recommend to prove following (known) geometrical identities
using this language in order to understand how useful it is:

\begin{lem}
The anticommutators can be computed by the following formulas:
 $$d\hat{X}+\hat{X}d=\nabla_{X}$$ $$ d\omega^* +\omega^*d
=(d\omega)$$ $$ \omega^*\hat{X}+\hat{X}
\omega^*=\omega(X)=(\omega_iX^i(x))$$ where $X$ is a vector field,
$\nabla_{X}$ is a Lie derivative of the differential forms along
the vector field $X$ and $\omega(X)$ is a value of 1-form on the
vector field.
\end{lem}

The proof easily follows  from the direct elementary calculation
with fermions. It is really much more convenient to calculate
everything with differential forms using the fermionic language.

So we have $$(d')^2=(d+\lambda\omega^* +\mu\hat{X})^2=
\lambda(d\omega)+\mu \nabla_{X}+\lambda\mu \omega(X)$$

\begin{lem} For the case $\lambda=0$  kernel of the operator $\nabla_{X}$ is
equal to the complex of differential forms invariant under the
one-parametric diffeomorphism group generated by the vector field
$X$. For $\mu=0$ corresponding homology group $H^*_{X,\mu=0}$ is
$Z$-graded and isomorphic to the homology of
$\nabla_{X}$-invariant De Rham Complex $H^*_{inv}$. In particular,
for the case of isometry on the compact manifold $M$ the homology
of this complex coincide with the standard $H^*(M)$.
\end{lem}

 In order to calculate the homology of perturbed
$Z_2$-graded complex with operator $d'=d+\mu\hat{X}$ we construct
a spectral sequence $E_l,d_l$ similar to the previous cases, using
the decomposition in the variable $\mu$. In the zero term we have
as before a $Z$-graded groups $$E_2=H^*_{X,\mu=0}=\sum_{i\geq 0}
H^i_{inv}, d_2(u)=\hat{X}(u)$$ for the representative closed forms
$u$. All higher differentials $d_l$ can be constructed easily:
they are based on the action of the operator $\hat{X}$ instead of
multiplication operator. This  operator is purely algebraic (non
differential). Last property allows us to call all higher
differentials the  analogs of Massey products.

\begin{lem} There exists a natural spectral sequence $E_l,d_l$ where
$$E_2=\sum_{i\geq 0} H^i_{inv}, d_2=\hat{X}:H^i_{inv}\rightarrow
H^{i-1}_{inv}$$ $$d_l:E_l^k\rightarrow E_l^{k-2l+3},l\geq 2$$ Here
$d_2$ is given by the contraction of the closed invariant form
with vector field $X$, i.e. $d_2(u)=\hat{X}(u)$. The higher
differentials $d_{l+1},l+1\geq 3$, are defined by the formula
$d_{l+1}(u)=\hat{X}(v)$ where $d_l(u)=dv$ for the $k$-cocycle $u$
presenting the corresponding element of the group
$E_{l+1}=\sum_{k\geq 0}E_{l+1}^k$. All $d_j$ with $j\leq l$
annihilate this element, and $d_{l+1}(u)$ is presented by the
$k-2l-1$-form in the class $d_{l+1}(u)\in E^{k-2l-1}_{l+1}$.
\end{lem}

{\it Remark: As before, we expect that in this case the group
$E_{\infty}$ for this spectral sequence is isomorphic to the
perturbed homology group of the operator $d+\mu\hat{X}$ for the
small $\mu\neq 0$. No problem to prove this fact for the case of
isometry group (i.e. if the closure of this one-parametric group
in the diffeomorphism group is compact). This case is similar to
the elliptic complexes. However, for the general case we will have
functional difficulties.}

Let us consider the special second differential
$$d_2=\hat{X}:H^1_{inv}\rightarrow H^0_{inv}=R$$ Following the old
works of Arnold and others of early 1960s, we consider an analog
of the ''rotation number'' for the vector field $X$ on the compact
manifold $M$. Starting from the fixed point $x$ belonging to the
topological closure of trajectory $x(t)$ of this dynamical system,
we define a limiting cycle $a\in H_1(M,R)$: consider a sequence
$t_k\rightarrow +\infty$ such that $x(t_k)\rightarrow x$. These
pieces of trajectory $[x(t_k),x(t_{k+1})]$ can be transformed into
the closed curves joining their ends by the small paths. We are
coming to the cycles $a_k\in H_1(M)$ made almost completely out of
pieces of the trajectory $x(t)$. There exists a limit
$$a_k/(t_{k+1}-t_k)\rightarrow a\in H_1(M), k\rightarrow +\infty$$
independent of trajectory.

\begin{lem}
Let $a=0$. Then the special differential
$d_2=\hat{X}:H^1_{inv}(M)\rightarrow R$ is equal to zero.
\end{lem}

Proof of this lemma easily follows from the fact that the value of
the invariant closed 1-form $u_1(X)$ on the trajectory is
constant. Therefore we have $d_2(u_1)=u_1(X)=(u_1,a)=0$. Lemma is
proved.

{\it Remark: If $a\neq 0$ and $X$ generates an isometry in some
Riemannian metric, we can prove an inverse statement.}

As a conclusion, let us mention that there are no natural
differential operators on the spaces of tensors independent on
Riemannian metric or other similar additional structures except
$d$ on the spaces of differential forms. However, there is a
natural metric independent class of purely algebraic operators
described above. All of them may be combined with $d$. As it was
demonstrated in the examples above, these exotic differentials and
complexes may lead to the interesting geometrical and analytical
objects.  One should not expect, however, that this sort of
quantities never appeared in topology before:

{\it Example: The Equivariant Homology}. Consider the operators
$d'=d+a_i\hat{X}_i$ for the polynomial generators $a_i$ in the
$Z_+$-graded algebra $R=R^+=Z(a_1,\dots,a_m)$ with 2-dimensional
polynomial generators $a_i$ This complex was invented by M.Atiyah
and others  for the definition of the ''equivariant homology''. A
lot of people studied them since that. Here $X_i$ are vector
fields (the infinitesimal isometries) on the manifold $M$
generating the action of compact commutative group $G=T^m$:
$[\nabla_{X_i},\nabla_{X_j}]=0$. The operator $d'$ acting on the
$G$-invariant $R$-valued forms can be naturally considered as a
perturbation of the ordinary $d$ in the subspace
$Ker(d')^2=\Lambda^*_{inv}\subset  \Lambda^*(M,R)$, so it is a
part of our general scheme: the spectral sequence for calculation
of these homology can be naturally considered as a by-product of
perturbations similar to the spectral sequence of Massey products.

\end{document}